\documentclass{article}
\usepackage{epsfig}

\usepackage{latexsym,amsfonts}

\usepackage{amsmath,amsthm}

\input xy
\xyoption{all}

\newcommand\frb{\mathfrak b}
\newcommand\frc{\mathfrak c}
\newcommand\frg{\mathfrak g}
\renewcommand\frm{\mathfrak m}
\newcommand\frn{\mathfrak n}
\newcommand\frp{\mathfrak p}

\newcommand\fru{\mathfrak u}

\newcommand\cA{\mathcal A}
\newcommand\cB{\mathcal B}
\newcommand\cF{\mathcal F}
\newcommand\cM{\mathcal M}
\newcommand\cN{\mathcal N}

\newcommand\cR{\mathcal R}

\newcommand\Aa{\mathbb A}
\newcommand\Ff{\mathbb F}
\newcommand\Zz{\mathbb Z}

\newcommand\lra{\longrightarrow}
\newcommand{\R}{T}

\DeclareMathOperator{\mcodim}{codim} %
\DeclareMathOperator{\mcd}{cd} %
\DeclareMathOperator{\mGL}{GL} %
\DeclareMathOperator{\mEnd}{End} %
\DeclareMathOperator{\mExt}{Ext} %
\DeclareMathOperator{\mext}{ext} %
\DeclareMathOperator{\mhom}{hom} %
\DeclareMathOperator{\mHom}{Hom} %
\DeclareMathOperator{\mIn}{in} %
\DeclareMathOperator{\mLie}{Lie} %
\DeclareMathOperator{\mm}{m} %
\DeclareMathOperator{\mSpec}{Spec} %
\DeclareMathOperator{\mStab}{Stab} %
\DeclareMathOperator{\msupp}{supp} %
\DeclareMathOperator{\mU}{U} %

\usepackage{amssymb}
\usepackage{amsfonts}

\numberwithin{equation}{section}

\newtheorem{Thm}{Theorem}[section]
\newtheorem{Lemma}[Thm]{Lemma}

\newtheorem{Prop}[Thm]{Proposition}

\theoremstyle{definition}
\newtheorem{Defn}[Thm]{Definition}
\newtheorem{Exmp}[Thm]{Example}
\newtheorem{Notn}[Thm]{Notation}

\theoremstyle{remark}
\newtheorem{Rem}[Thm]{Remark}

\begin{document}
\title{Prehomogeneous spaces for Borel subgroups of general linear groups}
\author{Simon M.~Goodwin and Lutz Hille}
\date{\today}

\maketitle

\begin{abstract}
Let $k$ be an algebraically closed field.  Let $B$ be the Borel
subgroup of $\mGL_n(k)$ consisting of nonsingular upper triangular
matrices. Let $\frb = \mLie B$ be the Lie algebra of upper
triangular $n \times n$ matrices and $\fru$ the Lie subalgebra of
$\frb$ consisting of strictly upper triangular matrices.  We
classify all Lie ideals $\frn$ of $\frb$, satisfying $\fru'
\subseteq \frn \subseteq \fru$, such that $B$ acts (by conjugation)
on $\frn$ with a dense orbit. Further, in case $B$ does not act with
a dense orbit, we give the minimal codimension of a $B$--orbit in
$\frn$.  This can be viewed as a first step towards the difficult
open problem of classifying of all ideals $\frn \subseteq \fru$ such
that $B$ acts on $\frn$ with a dense orbit.

The proofs of our main results require a translation into the
representation theory of a certain quasi-hereditary algebra
$\cA_{t,1}$.   In this setting we find the minimal dimension of
$\mExt^1_{\cA_{t,1}}(M,M)$ for a $\Delta$-good $\cA_{t,1}$--module
of certain fixed $\Delta$-dimension vectors.
\end{abstract}

\section{Introduction}

Let $k$ be an algebraically closed field. Let $P$ be a parabolic
subgroup of $\mGL_n(k)$, write $\frp$ for the Lie algebra of $P$ and
$\frp_u$ for the Lie algebra of the unipotent radical of $P$. The
group $P$ acts on $\frp_u$ via conjugation, and this induces an
action of $P$ on any Lie ideal $\frn$ of $\frp$ contained in
$\frp_u$. Thanks to a classical result of Richardson
(\cite{Richardson}) the parabolic group $P$ always acts on $\frp_u$
with a dense orbit. However, there is not always a dense $P$--orbit
in $\frn$, see \cite{HR3} for a counter example in case $\frn$ is
the derived subalgebra $\frp_u' = [\frp_u,\frp_u]$ of $\frp_u$.

Now consider the case $P=B$ is a Borel subgroup of $\mGL_n(k)$. We
write $U$ for the unipotent radical of $B$.  The Lie algebras of $B$
and $U$ are denoted by $\frb$ and $\fru$ respectively. We may take
$B$ to be the group of invertible upper triangular $n \times n$
matrices. Then $U$ consists of the upper unitriangular $n \times n$
matrices, $\frb$ is the Lie algebra of upper triangular $n \times n$
matrices, and $\fru \subseteq \frb$ consists of the strictly upper
triangular matrices. Let $\frn$ be an ideal of $\frb$ satisfying
$\fru' \subseteq \frn \subseteq \fru$. The principal result of this
paper gives a classification of all pairs $(B,\frn)$ such that
$\frn$ is a prehomogeneous space for $B$. We recall that $\frn$ is
called a {\em prehomogeneous space} for $B$ provided $B$ acts on
$\frn$ with a dense orbit.  This result can be viewed as a first
step towards classifying all ideals $\frn$ of $\frb$ contained in
$\fru$ such that $B$ acts on $\frn$ with a dense orbit.  This
classification is known to be a difficult problem, but we believe
that this paper contains many ideas and techniques for attacking
this problem.

Let $N$ be the normal subgroup of $B$ such that $U' \subseteq N
\subseteq U$, and let $\frn = \mLie N$. We note that if $x \in \frn$
is a representative of a dense $B$--orbit in $\frn$, then $1 + x$ is
a representative of a dense $B$--conjugacy class in $N$. Therefore,
our classification of when $B$ acts on $\frn$ with a dense orbit
gives a classification of when there is a dense $B$--conjugacy class
in $N$.

For $x \in \fru$, the orbit map $B \to B \cdot x$ is separable; this
follows from \cite[Prop.\ 6.7]{bor} and the fact that $C_B(x) = 1 +
\frc_\frb(x)$. Therefore, our classification of when $\frn$ is a
prehomogeneous space for $B$ leads to a classification of when there
exists $x \in \frn$ such that $[\frb,x] = \frn$. Moreover, since the
representative  $x \in \frn$ of the dense orbit (when it exists)
that we construct in Section \ref{sectrepresentative} is a matrix
whose entries are all 0 or 1, this classification remains valid if
we replace $k$ by an arbitrary field.

To formulate our main result we need to associate to the $B$--module
$\frn$ a sequence of positive natural numbers $\underline a =
(a_0,a_1,\ldots,a_r,a_{r+1})$.  The idea is that $\frn$ is
determined by which elementary matrices $e_{i,i+1}$ it contains. Let
$\{i_0,\dots,i_r\}$ be the set of $i$ such that $e_{i,i+1} \in
\frn$, where $i_j < i_{j+1}$ for each $j$.  Then we define
$\underline a$ by: $a_0 = i_0$; $a_j = i_j - i_{j-1}$, for $j =
1,\dots,r$; and $a_{r+1} = n - i_r$.  We refer the reader to Section
\ref{sectpreliminary} for an alternative definition of $\underline
a$ and an example.

We define $\mcd(B;\frn)$ to be the codimension of a $B$--orbit of
maximal dimension in $\frn$.  In our main theorem below, as well as
giving a classification of all instances when $B$ acts on $\frn$
with a dense orbit, we give the value of $\mcd(B;\frn)$ when $B$
does not act with a dense orbit.  In fact our proof constructs a
family of $B$--orbits of maximal dimension whose union is dense in
$\frn$.

\begin{Thm}\label{maintheorem}
Let $B$ be the Borel subgroup of $\mGL_n(k)$ consisting of upper
triangular matrices. Let $\frn$ be a Lie ideal of $\frb$ satisfying
$\fru' \subseteq \frn \subseteq \fru$ and associate the sequence
$\underline a$ to $\frn$ as above. Then: \\
{\em (1)} $B$ acts on $\frn$ with a dense orbit precisely when the
number of indices $i$ with $1 \le i \le r$ and $a_i$ even is at most
one.
\\
{\em (2)} If $B$ does not act on $\frn$ with a dense orbit, then the
minimal codimension of a $B$--orbit in $\frn$ is
$$
\mcd(B;\frn) = \sharp \{i \mid 1 \le i \le r, \, a_i \mbox{ is even
} \} - 1.
$$
\end{Thm}

The proof of the theorem uses methods both from algebraic group
theory, and from representation theory of algebras.  More
specifically we use algebraic group theoretic techniques to show
that certain ideals $\frn$ do not have a dense orbit, or the
codimension of a $B$--orbit of maximal dimension is at least the
claimed value.  In particular, we use results of the first author
from \cite{goo3} that give unique {\em minimal representatives} of
$B$--orbits in $\frn$.  This is necessary in order to show that
elements in the families that we construct are pairwise
non-conjugate.  For the converse we use representation theoretic
methods as explained in the next paragraph.

It was shown and intensively used in a series of papers that orbits
for parabolic group actions are related to modules over certain
quasi-hereditary algebras (see e.g.~\cite{HRTransgroups},
\cite{BHJAlg}, \cite{BHRR}, and \cite{BHRZ}). In our case we have to
deal with the algebra $\cA_{t,1}$. It was first considered in detail
in \cite{BHJAlg} (see Section \ref{sectpreliminary} for a definition
and relevant properties). For each element $x \in \frn$ there is a
corresponding $\cA_{t,1}$--module $M(x)$, so that elements in the
same $B$--orbit define isomorphic $\cA_{t,1}$--modules. For our
proofs a representation-theoretic characterization of a dense orbit
is crucial: an element $x$ in $\frn$ represents a dense orbit
precisely when $\mExt^1(M(x),M(x)) = 0$ (Theorem \ref{Thmext}). This
allows us to construct a dense orbit explicitly using certain
well-known standard $\cA_{t,1}$--modules, see Sections
\ref{sectdense}.  Further, we show in Theorem \ref{Thmext} that the
codimension of a $B$--orbit in $\frn$ is equal to the dimension of
$\mExt^1(M(x),M(x))$.   This result is of importance for our
construction of a family of $B$--orbits whose union is dense in
$\frn$.



\medskip

As hinted at above there has been a lot of recent interest in the
action of a parabolic subgroup $P$ of a reductive algebraic group
$G$ on the Lie algebra $\frp_u$ of its unipotent radical.  For
example, there is a classification of all instances when $P$ acts on
$\frp_u$ with finitely many orbits, see \cite{HRTransgroups} and
\cite{jurroh}. Further, for $G$ not of type $E_7$ or $E_8$, there is
a classification of all instances when $P$ acts on higher terms
$\frp_u^{(l)}$ of the descending central series of $\frp_u$ with a
finite number of orbits, see \cite{BHJAlg}, \cite{BHR} and
\cite{gooroh2}.

There has also been much interest in the question of when
$\frp_u^{(l)}$ is a prehomogeneous space for $P$. This was first
investigated by G.~R\"ohrle and the second author in \cite{HR3}. It
was considered further by G.~R\"ohrle and the first author in
\cite{gooroh1} and \cite{gooroh2}.  In particular, the results of
these papers give a classification of all instances when a Borel
subgroup $B$ acts of $\fru^{(l)}$ with a dense orbit. Further, in
\cite{goo2} the first author describes a computer program that
determines when a Borel subgroup $B$ acts on any $B$--submodule of
$\fru$ with a dense orbit.

There has been other interest in the adjoint orbits of a Borel
subgroup $B$ of a reductive algebraic group on $\fru$.  The problem
of determining the $B$--orbits (and $U$--orbits) in $\fru$ was
addressed by H.\ B\"urgstein and W.H.\ Hesselink in \cite{burhes}.
They were motivated by the problem of describing the component
configuration of the {\em Springer fibre} $\cB_x = \{B \in \cB : x
\in \mLie B \}$ where $\cB$ denotes the variety of Borel subgroups
of $G$ and $x \in \frg$ is nilpotent, see \S 1.4 of {\em loc.\ cit.\
} The first author studied the $B$--orbits and $U$--orbits in $\fru$
in \cite{goo3}.  In particular, in {\em loc.\ cit.\ }it is shown
that any $U$--orbit contains a unique so-called {\em minimal
representative}.

\smallskip

In case $G$ is a reductive algebraic group defined and split over
the finite field $\Ff_q$ of $q$ elements ($q = p^s$, where $p$ is
prime), there has been much interest in the conjugacy classes of a
Sylow $p$-subgroup of the Chevalley group $G(q)$.  If $B$ is a Borel
subgroup of $G$ that is defined over $\Ff_q$, then $U(q)$ is a Sylow
$p$-subgroup of $G$. For the case $G = \mGL_n(k)$, we refer the
reader to \cite{hig}, \cite{rob}, \cite{tho} and \cite{verarr}; for
arbitrary $G$, see \cite{goo3}.

In a final section of this paper we give an interpretation of
Theorem \ref{maintheorem} in the corresponding finite group setting.
In Theorem \ref{thmfinite} we give the maximal size of a
$B(q)$--conjugacy class in a normal subgroup $N(q)$ of $B(q)$ such
that $U'(q) \subseteq N(q) \subseteq U(q)$.  The proof of this
theorem requires results from \cite{goo3}.

\medskip

We now outline the structure of this paper.  After giving some
notation below we cover the requisite preliminaries in Section
\ref{sectpreliminary}.  We give the basic notation we use for the
action of $B$ on $\frn$ and then we discuss the algebras
$\cA_{t,l}$.  In particular, we give some results for calculating
$\mHom$-groups and $\mExt^1$-groups for {\em standard}
$\cA_{t,1}$--modules that we require in the sequel.  We explain the
relationship between $P$--orbits in $\frp_u$ with $\Delta$--good
modules for the algebra in $\cA_{t,l}$ and prove Theorem
\ref{Thmext} in Section \ref{sectext}.  Next in Section
\ref{sectdense} we prove one direction of Theorem
\ref{maintheorem}(1), by constructing certain $\Delta$--good
$\cA_{t,1}$--modules without self-extensions.  We give a
construction of a representative of the dense $B$--orbit in $\frn$
in Section \ref{sectrepresentative} and prove that it is minimal in
the sense of \cite{goo3}; the required terminology and results from
{\em loc.\ cit.\ }are recalled.  These representatives are required
for our proof, in Section \ref{sectnotdense}, that there is not a
dense $B$--orbit in $\frn$ in the cases stated in Theorem
\ref{maintheorem}(1).   In Section \ref{sectparameters}, we prove
Theorem \ref{maintheorem}(2).  This is achieved by first
constructing certain $\Delta$-good $\cA_{t,1}$--modules $M$ for
which we can give an upper bound on the dimension of $\mExt^1(M,M)$.
Then we construct a family $\cF$ of elements of $\frn$ that are
pairwise non-conjugate and for which we can bound the dimension of
$\mExt^1(M(x),M(x))$, where $x$ is a geometric fibre of $\cF$.
Dimension arguments allow us to show that the family of $B$--orbits
given by $\cF$ has dense union in $\frn$ completing the proof of
Theorem \ref{maintheorem}(2).  The final section is dedicated to
proving Theorem \ref{thmfinite}.

\subsection*{Basic Notation and references}
With $\subset$ we denote the strict inclusion and with $\sharp$ the
number of elements in a finite set.  Given $x \in \mathbb{R}$ we
write $\lfloor x \rfloor$ for the greatest integer less than $x$.

We always work over an algebraically closed field $k$. We write
$k^\times$ for the nonzero elements of $k$. All algebras are
$k$--algebras with unit and modules are finitely generated.
We denote by $\otimes$ the tensor product over $k$. If we consider
families of modules or families of elements we work over finitely
generated commutative $k$--algebras.  We call such an algebra an
{\em affine $k$--algebra}; it is, in general, not reduced.

An {\em affine family} $\cM$ of $\cA$--modules is just an
$\cA$--$R$--bimodule, locally free over $R$, where in our situation
$R$ is an affine $k$--algebra and $\cA$ is a finite dimensional
$k$--algebra.  Given a maximal ideal $m$ of $R$, we have an
$\cA$--module $\cM(m) = \cM \otimes_R R/m$.  In this way $\cM$
becomes a family of modules over the maximal spectrum $\mSpec_{\max}
R$ of $R$.

An {\em affine family $\cF$ of elements in $\frp_u$} is just an
element in $\frp_u \otimes R$. If $m \subset R$ is a maximal ideal,
then the projection $R \to R/m \cong k$ induces a morphism $\frp_u
\otimes R \to \frp_u \otimes R/m = \frp_u$. Given $\cF \in\frp_u
\otimes R$ we write $\cF(m) = \cF \otimes_R R/m$ for its image in
$\frp_u$.   In this way $\cF$ becomes a family of elements over
$\mSpec_{\max} R$.

If we write $\mHom$ or $\mExt^1$ we always consider morphisms or
extensions of modules for the algebra $\cA_{t,l}$. We write $\mhom$
to denote the dimension of a $\mHom$-group and $\mext^1$ to denote
the dimension of an $\mExt^1$-group.

For some background on path algebras and quivers we refer to
\cite{GabrielRoiter} and \cite{Ringel1}. Basic results on
quasi-hereditary algebras can be found in \cite{DlabRingel}. The
quasi-hereditary algebras $\cA_{t,l}$ together with certain
particular modules were considered in \cite{BHJAlg}.   First results
concerning actions with dense orbits and tilting modules in our
setting can be found in \cite{BHRR} and \cite{HR3}.  As a reference
for the theory of algebraic groups we refer the reader to
\cite{bor}.

\section{Preliminary results}\label{sectpreliminary}

Let $\{0\} = V_0 \subseteq V_1 \subseteq V_2 \subseteq \ldots
V_{t-1} \subseteq V_t = V$ be a flag in the finite dimensional
vector space $V$ over $k$. We define the dimension vector of this
flag to be $d=(\dim V_1, \dim V_2 - \dim V_1,\ldots,\dim V_t - \dim
V_{t-1})$. We also define the dimension vector $\Sigma d = (\dim
V_1,\ldots,\dim V_t)$. So we have that $(\Sigma d)_i = \sum_{j=1}^i
d_j$. In this paper we allow also some $d_i$ equal zero, so the flag
is not proper in general.

Let $P = P(d) \subseteq \mGL(V)$ be the stabilizer of the flag and
$\frp = \frp(d)$ the Lie algebra of $P$.  The Lie algebra of the
unipotent radical of $P$ is given by $\frp_u = \frp_u(d) = \{ f \in
\mEnd(V) \mid f(V_i) \subseteq V_{i-1} \}$.  We define the ideal
$\frp_u^{(1)} = \frp_u^{(1)}(d)$ of $\frp$ by $\frp_u^{(1)} = \{ f
\in \mEnd(V) \mid f(V_i) \subseteq V_{i-2} \}$. We note that the
group $P$ and the Lie algebra $\frp_u$ do not depend on how many
zeros we include in $d$. However, the ideal $\frp_u^{(1)}$ differs
if we insert zeros in $d$.

In this note we are interested in dimension vectors $d$ with $d_i =
0$ or $d_i = 1$, for all $i$; such a dimension vector is called {\em
thin}. For thin dimension vectors we have that $P = B$ is a Borel
subgroup of $\mGL(V)$ and we write $\fru = \frp_u$ and $\frn =
\frp_u^{(1)}$.  It is straightforward to show that any ideal $\frn$
of $\frb$ satisfying $\fru' \subseteq \frn \subseteq \fru$, can be
obtained in this manner for some $d$, see \cite[Satz
1.4.2]{Hillehabil}.

In case $d$ is a thin dimension vector, we can make the following
assumptions on $d$ without affecting $B$ and $\frn$: $d_1=1$;
$d_t=1$; and if $d_i=0$, then $d_{i+1}=1$. So we make these
assumptions for the rest of the paper without loss of generality. In
the definition below we associate to a thin dimension vector, the
sequence $\underline a$ counting the lengths of the strings of 1 in
$d$; we note that this definition of $\underline a$ is equivalent to
that given in the introduction.

\begin{Defn} \label{defna(d)}
Let $d$ be a thin dimension vector.  We define the sequence
$\underline a = \underline a(d) = (a_0,a_1,\ldots,a_{r},a_{r+1})$ of
positive integers by setting $a_i$ to be the length of the $i$th
string of ones in $d$.  That is $a_0+1$ is the position of the first
entry 0 in $d$, $a_0 + a_1 + 2$ is the position of the second entry
0 in $d$, and so on.  The nonnegative integer $e(d)$ is defined by
$$
e(d) = \sharp \{i \mid 1 \le i \le r, \, a_i \mbox{ is even } \}.
$$
\end{Defn}

We note that conversely each sequence $\underline a =
(a_0,a_1,\ldots,a_r,a_{r+1})$ of positive integers defines a thin
dimension vector $d$
starting with a
sequence of $a_0$ ones followed by a zero, a sequence of $a_1$ ones
followed by zero until the last sequence of ones of length
$a_{r+1}$.

\begin{Exmp} We consider the case where $t=17$, $\dim V = 13$ and $d =
(1,1,0,1,1,1,0,1,1,0,1,1,1,1,1,0,1)$.  Then $r=3$ and $\underline a
= (2,3,2,5,1)$. By Theorem \ref{maintheorem} the group $B \subset
\mGL(V)$ acts on $\frn$ with a dense orbit.

\begin{figure} [h!tb]

\begin{center}

\includegraphics[width=5cm]{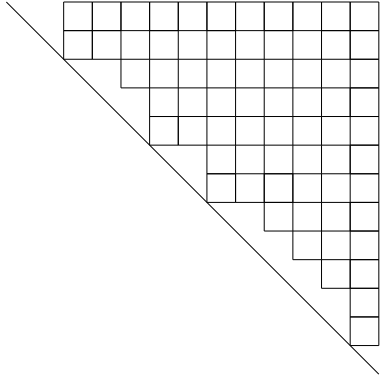}

\caption{The ideal $\frn $  for  $\underline a =
  (2,3,2,5,1)$} \label{f:1}

\end{center}

\end{figure}

\end{Exmp}

\begin{Notn} \label{notn}
We give some notational conventions that we use throughout the
sequel.

Let $d = (d_1,\dots,d_t)$ be a thin dimension vector and let $n =
\sum_{i=1}^t d_i$; as mentioned above we always assume that there
are no consecutive 0 entries in $d$ and that $d_1 = d_t = 1$.  The
support of $d$ is defined by $\msupp d = \{i \mid d_i = 1\}$.

We define the standard flag $F(d) : V_0 \subseteq \dots \subseteq
V_t$ corresponding to $d$ by setting $V_i = k^{(\Sigma d)_i}$. In
this way $d$ defines a Borel subgroup $B = P(d) = \mStab(F(d))$ of
$\mGL_n(k)$. Further, $d$ determines the Lie ideal $\frn =
\frp_u^{(1)}$ of $\frb$.  In the sequel we often refer to $B$ and
$\frn$ as the Borel subgroup and Lie ideal corresponding to $d$.

In order to make some constructions and proofs in the sequel easier
to understand we relabel the standard basis $\{e_1,\dots,e_n\}$ of
$k^n$ as follows: we define the function $\gamma: \msupp d \to
\{1,\dots,n\}$ by requiring $j$ is the $\gamma(j)$th lowest element
in $\msupp d$, and setting $f_j = e_{\gamma(j)}$.  For $i,j \in
\msupp d$ we let $f_{ij}$ be the elementary matrix that sends $f_i$
to $f_j$ and all other $f_{i'}$ to zero.
\end{Notn}

In the remainder of this section we discuss an algebra $\cA_{t,l}$.
In the next section we explain its connection with orbits of
parabolic subgroups that was established in \cite{BHparact}.  The
algebras $\cA_{t,l}$ are shown to be quasi-hereditary in {\em loc.\
cit.\ }as is reflected in the some of the terminology given below.

For positive integers $t$ and $l$ with $l < t$, we define a quiver
$Q_{t,l}$ as follows: the vertex set is $\{1,2,\ldots,t\}$; and
there are arrows $\alpha_i:i \lra i+1$ for $i= 1,\dots,t-1$ and
$\beta_j: j+l+1 \lra j$ for $j = 1,\dots,t-l-1$. The algebra
$\cA_{t,l}$ is defined to be the quotient of the path algebra
$kQ_{t,l}$ by the relations:
$$
\beta_1\alpha_{l+1} = 0  \text{ and } \beta_{i+1}\alpha_{i+l+1} =
\alpha_{i} \beta_{i} \text{, for } i=1,\ldots,t-l-2.
$$
We identify $\cA_{t,l}$--modules with the corresponding
representations of the quiver $Q_{t,l}$.  So we write $M_i$ for the
vector subspace of $M$ corresponding to the vertex $i$, and by an
abuse of notation given an arrow $\gamma : i \to j$ in $Q_{t,l}$ we
also write $\gamma : M_i \to M_j$ for the corresponding linear map.

\begin{figure} [h!tb]

\begin{center}

\centerline{ \includegraphics[width=10cm]{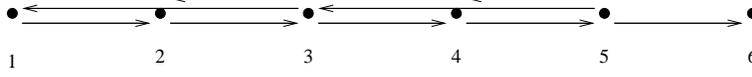}}

\caption{The quiver $Q_{6,1}$}

\end{center}

\end{figure}

The algebra $\cA_{t,l}$ was introduced in \cite{BHparact}, where it
was observed to be quasi-hereditary. It can also be realized as a
subalgebra of the endomorphism algebra of $\bigoplus_{i=0}^t
k[T]/T^i$ as a $k[T]$--module.
For the remainder of this section we specialize to the case $l = 1$,
which is the case in which we are interested in this paper; much of
the discussion below generalizes easily to the case $l \ge 2$.

We need several particular $\cA_{t,1}$--modules. Denote by
$\Delta(i)$ the representation with $\Delta(i)_j = k$ for $j \geq i$
and $\Delta(i)_j = 0$ otherwise, where all the arrows $\beta$ act
with zero and the arrows $\alpha $ with identity, whenever possible.
The projective cover of $\Delta(i)$ is denoted by $P(i)$.  Let
$e_1,\dots,e_m$ be the standard basis for $k^m$, for $m \in \Zz_{\ge
0}$. Then the projective module $P(i)$ is given by: $P(i)_j = k^l$,
where $l = \lfloor (j+2-i)/2 \rfloor$ if $j \ge i$ and $P(i)_j = 0$
for $j < i$; $\alpha_i$ is the natural inclusion; and $\beta_j(e_h)
= e_{h-1}$ if $h \ge 2$ and $\beta_j(e_1) = 0$.

Let $M$ be an $\cA_{t,1}$--module. We denote the dimension vector of
$M$ by $\dim M = (\dim M_1,\dots,\dim M_t)$. A module $M$ admitting
a filtration with each factor of the form $\Delta(i)$ for some $i$
is called $\Delta$-{\em good}. For such a module $M$ we  define its
$\Delta$-dimension vector $\dim_{\Delta} M = (\dim M_1, \dim M_2 -
\dim M_1,\ldots,\dim M_t - \dim M_{t-1})$. This vector counts the
number of filtration factors isomorphic to $\Delta(i)$. For example
$\dim_{\Delta}\Delta(i) = (0,\ldots,0,1,0,\ldots,0)$ and
$\dim_{\Delta} P(i) = (\dots,1,0,1,0,1,0,0,0,\ldots,0,0)$. We also
define the $\Delta$-support of a $\Delta$-good module $M$ as the set
of all $i$ with $(\dim_{\Delta} M)_i \not= 0$.

Let $J$ be a subset of $\{1,2,\ldots,t\}$ such that for each $j \in
J$ we have that $j+1 \notin J$; we call such $J$ a {\em standard
subset} of $\{1,\dots,t\}$. Let $J = \{j_1,\dots,j_r\}$ be standard
subset of $\{1,\dots,t\}$, where $j_i < j_{i+1}$ for each $i$. Then
there exists a unique indecomposable $\Delta$-good module
$\Delta(J)$ with a thin $\Delta$-dimension vector and
$\Delta$-support $J$ called a {\em standard module}; $\Delta(J)$ is
defined by: $\Delta(J)_i = k^l$, where $l$ satisfies $j_l \le i <
j_{l+1}$ (by convention we set $j_0 = 0$ and $j_{r+1} = t+1$);
$\alpha_i$ is the natural inclusion; and $\beta_j(e_h) = e_{h-1}$ if
$h \ge 2$ and $\beta_j(e_1) = 0$. So if $J$ is of the form $J =
\{\dots,i-4,i-2,i\}$, then we have $\Delta(J) = P(i)$, and if $J =
\{i\}$, then $\Delta(J) = \Delta(i)$.

In the following theorem we collect some properties of $\Delta$-good
modules for $\cA_{t,1}$.  Parts (1) and (2) of Theorem
\ref{ThmDeltagoodproperties} can be found in \cite{BHJAlg}. The last
claim (3) is a generalization of \cite[Prop.\ 1]{BHRR}; it is also a
consequence of Lemma \ref{lemdimext} below.

\begin{Thm}\label{ThmDeltagoodproperties}
{\em (1)} For an $\cA_{t,l}$--module $M$ the following conditions
are equivalent

{\em (a)} $M$ is $\Delta$-good;

{\em (b)} the linear maps $\alpha_i: M_i \to M_{i+1}$ are injective
for all $i=1,\ldots,t-l-1$;

{\em (c)} the projective dimension of $M$ is at most one; and

{\em (d)} the socle of $M$ is isomorphic to a direct sum of copies
of $\Delta(t)$.
\\
{\em (2)} The projective modules $P(t-l),P(t-l-1),\dots,P(t)$ are
also injective. \\
{\em (3)} Let $l = 1$  and let $J$ be a standard subset of
$\{1,\dots,t\}$. Then we have $\mExt^1(\Delta(J),\Delta(J)) = 0$.
\end{Thm}

In the sequel we need to know the dimensions of $\mHom$-groups
between standard modules.  These dimensions are given in the
following lemma.

\begin{Lemma} \label{lemdimhom}
Let $J = \{j_1,\dots,j_r\}$ and $K = \{k_1,\dots,k_s\}$ be standard
subsets of $\{1,\dots,t\}$, where $j_i < j_{i+1}$, $k_i < k_{i+1}$
for each $i$. Then $\hom(\Delta(J),\Delta(K))$ is the greatest $m
\le s$ such that $r \ge m$ and $j_{r-i} \ge k_{m-i}$ for $i =
0,1,\dots,m-1$.
\end{Lemma}

\begin{proof}
We have $\Delta(J)_t = k^r$ and $\Delta(K)_t = k^s$.  We see that
any homomorphism $\phi : \Delta(J) \to \Delta(K)$ is determined by
the map $\phi_t : \Delta(J)_t \to \Delta(K)_t$ and in turn this map
is determined by $\phi_t(e_r) = \sum_{i=1}^m b_i e_i \in
\Delta(K)_t$, where $m$ is such that $b_m \ne 0$.  For $\phi$ to be
a homomorphism of $\cA_{t,1}$--modules, $\phi_t$ must restrict to
the maps $\phi_i : \Delta(J)_i \to \Delta(K)_i$ for each $i$.  It is
easy to check that this holds if and only if $m$ must satisfy the
conditions of the lemma.  It follows that
$\hom(\Delta(J),\Delta(K))$ is the greatest such $m$.
\end{proof}

The proof of the lemma above allows us to define certain morphisms
between standard modules.

\begin{Defn} \label{defmor}
Let $J = \{j_1,\dots,j_r\}$ and $K = \{k_1,\dots,k_s\}$ be standard
subsets of $\{1,\dots,t\}$, where $j_i < j_{i+1}$, $k_i < k_{i+1}$
for each $i$.  Let $m \le s$ be such that $r \ge m$ and $j_{r-i} \ge
k_{m-i}$ for $i = 0,1,\dots,m-1$.  We define $\phi^{k_m} : \Delta(J)
\to \Delta(K)$ to be the unique morphism satisfying
$(\phi^{k_m})_t(e_r) = e_m$.
\end{Defn}

We give an example of the notation below.

\begin{Exmp}
Let $t = 7$, and let $J$ and $K$ be the standard subsets $\{1,3,7\}$
and $\{2,4,6\}$ respectively.  A morphism $\phi : \Delta(J) \to
\Delta(K)$ is determined by $\phi_7(e_3)$.  We cannot have $\phi_7(e_3)
= e_3$, because then we have $\phi_3(e_2) = e_2$, which is not
possible, because $e_2 \notin \Delta(K)_3$. The
morphism $\phi^4$ is given by $\phi^4_7(e_3) = e_2$, $\phi^4_7(e_2) =
e_1$ and $\phi^4_7(e_1) = 0$; and the morphism $\phi^2$ is given by
$\phi^2_7(e_3) = e_1$ and $\phi^2_7(e_2) = \phi^2_7(e_1) = 0$.
\end{Exmp}

We also require the following lemma regarding the vanishing of
certain $\mExt^1$-groups between standard modules.

\begin{Lemma} \label{lemdimext}
Let $J$ and $K$ be standard subsets of $\{1,\dots,t\}$.  Assume that
if $i \in K$, then $i-1 \in J$ or $i \in J$.  Then
$\mExt^1(\Delta(K),\Delta(J)) = 0$.
\end{Lemma}

\begin{proof}
Since $\Delta(K)$ is filtered by the standard modules $\Delta(i)$
for $i \in K$, it is sufficient to prove the lemma for $K = \{i\}$.
In this case we have a projective resolution
$$
0 \lra P(i-2) \overset{\psi}{\lra} P(i) \lra 0
$$
of $\Delta(i)$, where $\psi = \phi^{i-2}$.  This gives rise to the
long exact sequence
$$
0 \lra \mHom(\Delta(i),\Delta(J)) \lra \mHom(P(i),\Delta(J)) \lra
$$
$$
\qquad \qquad \lra \mHom(P(i-2),\Delta(J)) \lra
\mExt^1(\Delta(i),\Delta(J)) \lra 0.
$$
Now using Lemma \ref{lemdimhom}, and recalling that projective
modules are certain standard modules, we have
$$
\mhom(\Delta(i),\Delta(J)) = 1 \text{, } \mhom(P(i),\Delta(J)) = l
\text{ and } \mhom(P(i-2),\Delta(J)) = l-1,
$$
where $l = \sharp(J \cap \{1,\dots,i\})$.  We deduce that
$\mExt^1(\Delta(i),\Delta(J)) = 0$ as required.
\end{proof}

\section{The codimension of an $P$--orbit and $\mExt^1$}
\label{sectext}


Let $d = (d_1,\dots,d_t)$ be a dimension vector with $\sum_{i=1}^t
d_i = n$; we do not assume that $d$ is thin.  Then we define the
standard flag $F(d) : V_0 \subseteq \dots \subseteq V_t$ in $k^n$,
by $V_i = k^{(\Sigma d)_i}$.  The parabolic subgroup $P = P(d)$ of
$\mGL_n(k)$ is defined as in the previous section.  For each $l \ge
0$ we define the ideal $\frp_u^{(l)}(d) = \{ f \in \mEnd(V) \mid
f(V_i) \subseteq V_{i-l-1} \}$.

In what follows we need to compare $P$--orbits and $\Delta$-good
modules. The main theorem in this section (Theorem \ref{Thmext})
compares the codimension of a $P(d)$-orbit in $\frp_u^{(l)}$ with
the group of self-extensions of a corresponding $\Delta$-good
$\cA_{t,l}$--module $M(x)$.

We define $\cF_{t,l}(\Delta)$ to be the category of $\Delta$-good
finite dimensional modules for $\cA_{t,l}$ (the parameters $t$ and
$l$ are always fixed). Further we define $\cR(\cF_{t,l}(\Delta),d)$
to be the space of all $\Delta$-good representations of $\cA_{t,l}$
with $\Delta$-dimension vector $d$
$$
\cR(\cF_{t,l}(\Delta); d) = \{ (\phi_i,\psi_j) \in
\bigoplus_{i=1}^{t-1} \mHom(V_i,V_{i+1}) \oplus
\bigoplus_{j=1}^{t-1-l} \mHom(V_{j+l+1},V_{j})\mid \hspace{1cm} ~
$$
$$
~ \hspace{3cm} \phi_i \mbox{ is
  injective } \psi_1\phi_{l+1} = 0, \psi_{i+1}\phi_{i+l+1} = \phi_{i}
\psi_{i} \}.
$$
It is a scheme over $k$ and the group $\mGL(\Sigma d) =
\prod_{i=1}^t \mGL((\Sigma d)_i)$ acts via conjugation. Note that
the dimension vector of a $\Delta$-good module with
$\Delta$-dimension vector $d$ is just $\Sigma d$. The orbits of the
underlying affine variety under this action are in natural bijection
with the isomorphism classes of  $\Delta$-good $\cA_{t,l}$--modules
with $\Delta$-dimension vector $d$.

Below we outline the correspondence of $P(d)$--orbits in
$\frp_u^{(l)}$ and isoclasses of $\Delta$-good modules with
$\Delta$-dimension vector $d$ that was established in
\cite{BHparact}.

Let $x$ be an element in $\frp_u^{(l)}(d)$ for some dimension vector
$d$. We can associate to $x$ a $\Delta$-good module $M(x)$ over
$\cA_{t,l}$ of dimension vector $\Sigma d$ -- equivalently an
element of $\cR(\cF_{t,l}(\Delta); d)$. Note that $x$ maps $V_i$ to
$V_{i-1-l}$. So we define $\phi_i$ to be the natural inclusion $V_i
\lra V_{i+1}$ and $\psi_i$ the restriction of $x$ to $V_{i+l+1}$. It
is obvious that the maps $\phi_i$ and $\psi_j$ satisfy the necessary
conditions.

Conversely, given $(\phi_i,\psi_i) \in \cR(\cF_{t,l}(\Delta); d)$.
The map $\psi_{t-l-1}:V_t \lra V_{t-l-1}$ is an element in $\tilde
\frp_u^{(l)}(d)$ for the stabilizer $\tilde P$ of the flag in $V =
V_t$ defined by the images of the injections $\phi_i$. So the
various $\Delta$-good modules $M$ define elements $x \in \tilde
\frp_u^{(l)}(d)$ for possibly different groups $\tilde P$.

\medskip

For an affine $k$--algebra $R$, we consider the set of $R$--valued
points of the scheme $\cR(\cF_{t,l}(\Delta); d)$, which is denoted
by $\cR(\cF_{t,l}(\Delta); d)(R)$.  It consists of all pairs
$$
(\phi_i,\psi_i) \in \left( \bigoplus_{i=1}^{t-1} \mHom(V_i,V_{i+1})
\otimes R \right) \oplus \left( \bigoplus_{j=1}^{t-1-l}
\mHom(V_{j+l+1},V_{j}) \otimes R \right)
$$
satisfying the following two conditions:
\begin{enumerate}
\item[(1)] the maps $\phi_i$ are split injective; and
\item[(2)] $\psi_1\phi_{l+1} = 0$ and $\psi_{i+1}\phi_{i+l+1} =
\phi_{i} \psi_{i}$.
\end{enumerate}

A family (over $R$) of $\Delta$-good modules is an
$\cA_{t,l}$--$R$--bimodule $\cM$, so that the restriction of $\cM$
to the directed quiver of type $\Aa_t$ (defined by the arrows
$\alpha_i: i \lra i+1$, see Figure 2) consists of locally direct
maps of locally free $R$--modules. This condition can also be read
as follows: for each maximal ideal $m$ in $R$ there exists a
localization $R_f$, with $f$ not in $m$, and an element $g \in
\mGL(\Sigma d)(R_f) = \prod_{i=1}^t \mGL((\Sigma d)_i)(R_f)$ so that
$g\cdot(\cM(\alpha_i))$ is the standard flag $F(d)$ (considered over
$R_f$). So after passing to a localization $R_f$ of $R$ and choosing
a basis of the free $R_f$--module $\cM$, we can assume $\cM \otimes
R_f$ is an element in $\cR(\cF_{t,l}(\Delta); d)(R_f)$. If, in
particular, $R$ is a local ring (for example the ring of dual
numbers), then isomorphism classes of families are in one-to-one
correspondence with orbits in $\cR(\cF_{t,l}(\Delta); d)(R)$ under
the action of $\mGL(\Sigma d)(R)$.

Let $x$ be an element in $\frp_u^{(l)}(d)$ and consider a family
$\overline x$ in $\frp_u^{(l)}(d) \otimes
k[\varepsilon]/\varepsilon^2$ with $\overline x
\otimes_{k[\varepsilon]/\varepsilon^2} k = x(\varepsilon) = x$, then
$\overline x$ is just a tangent vector in $x$, so an element in
$T_x(\frp_u^{(l)}(d))$. Then the group $\mGL(\Sigma
d)(k[\varepsilon]/\varepsilon^2)$ acts on the space of all families
$\cF$ of elements in $\frp_u^{(l)}(d)$ over the ring
$k[\varepsilon]/\varepsilon^2$. The orbit under this action
intersected with $T_x(\frp_u^{(l)}(d))$ can be naturally be
identified with the quotient $T_x(\frp_u^{(l)}(d))/T_x(P \cdot x)$.
We can consider the analogous families of $\Delta$-good
$\cA_{t,l}$--modules. Then the main step in the proof of Theorem
\ref{Thmext} consists of the identification of the respective
quotients of the tangent spaces.

The following result shows in particular that for a family of
$\Delta$-good $\cA_{t,l}$--modules we can at least on a sufficiently
small Zariski open  subset construct a family of elements in
$\frp_u^{(l)}(d)$. Moreover, for $x \in \frp_u^{(l)}(d)$ the
codimension of the $P$--orbit of $x$ in $\frp_u^{(l)}(d)$ coincides
with the dimension of the extension group $\mExt^1(M(x),M(x))$. This
result was first proven in \cite{Hillehabil}. The Artin--Voigt Lemma
(see \cite[\S 1]{Gabrielopen}) tells us that $\dim
\mExt^1(M(x),M(x))$ coincides with the dimension of the quotient of
the scheme-theoretic tangent space of the representation space of
all $\cA_{t,l}$--modules in $M(x)$ by the tangent space along the
$\mGL(\Sigma d)$--orbit. Since $\cR(\cF_{t,l}(\Delta); d)$ is an
open subscheme we can, for the $\Delta$-good module $M(x)$, also
consider the dimension of $T_{M(x)}(\cR(\cF_{t,l}(\Delta));
d)/T_{M(x)}(\mGL(\Sigma d) \cdot M(x))$.

\begin{Thm}\label{Thmext}
Let $x$ be an element in $\frp_u^{(l)}(d)$ and $M(x)$ the
corresponding $\Delta$-good $\cA_{t,l}$--module. Then the
codimension of the orbit $P\cdot x$ in $\frp_u^{(l)}(d)$ coincides
with the dimension of $\mExt^1(M(x),M(x))$. In particular, the
element $x$ is a representative of a dense orbit precisely when
$\mExt^1(M(x),M(x))= 0$.
\end{Thm}

\begin{proof} We consider an affine family $\cF$ of elements in
$\frp_u^{(l)}(d)$: this is just an element in
$\frp_u^{(l)}(d)\otimes R$ for an affine $k$--algebra $R$. To prove
the theorem we only require the case where $R$ is the ring of dual
numbers $k[\varepsilon]/ \varepsilon^2$;
however the arguments apply for any affine $k$--algebra $R$. Such an
affine family $\cF$ defines an $\cA_{t,l}$--module $\cM$, free over
$R$, as follows: for $\phi_i$ we can take the inclusion $V_i \to
V_{i+1}$ (and extend it to $V_i \otimes R \subseteq V_{i+1} \otimes
R$) and we define $\psi_i = \cF|_{V_{i+1+l} \otimes_{k} R} \in
\mHom(V_{i+1+l}, V_i) \otimes R$. So it is an element in $(\bigoplus
\mHom(V_i,V_{i+1}) \oplus \bigoplus \mHom(V_{j+l+1},V_{j})) \otimes
R$ satisfying the same conditions (relations and injectivity) like a
$\Delta$-good representation of $\cA_{t,l}$. Moreover, the maps
$\phi_i$ are obviously locally direct (they are in fact split over
$R$). In this way each affine family of elements in
$\frp_u^{(l)}(d)$ defines a corresponding family of
$\cA_{t,l}$--modules.

Conversely, we also need to show that for each affine family $\cM$
of $\cA_{t,l}$--modules, we can cover $\cM$ by families $\cN$ for
which there is a corresponding family $\cF$ of elements in
$\frp_u^{(l)}(d)$.  We cannot expect to get a corresponding family
of elements in $\frp_u^{(l)}(d)$ for all affine families of
$\cA_{t,l}$--modules.

In order to see this we consider the affine space
$$
\cR(\Aa_t; d) = \left\{ (\phi_i) \in \bigoplus_{i=1}^{t-1}
\mHom(V_i,V_{i+1}) \,\, \vline \,\, \dim V_i = (\Sigma d)_i \right\}
$$
of representations of a directed quiver of type $\Aa_t$. It contains
the open subscheme
$$
\cR(\Aa_t; d)^{\mathrm{inj}} = \left\{ (\phi_i) \in
\bigoplus_{i=1}^{t-1} \mHom(V_i,V_{i+1}) \,\, \vline \,\, \dim V_i =
(\Sigma d)_i, \phi_i \mbox{ is
  injective } \right\}.
$$
and there is a natural projection morphism
$$
\pi: \cR(\cF_{t,l}(\Delta); d) \lra \cR(\Aa_t; d)^{\mathrm{inj}}, \quad
(\phi_i,\psi_j) \mapsto (\phi_i).
$$
Then the $R$--valued points $\cR(\Aa_t; d)^{\mathrm{inj}}(R)$ consist
of split injective maps $\phi = (\phi_i)$. In particular,
we can find an element $ g \in \mGL(\Sigma d)(R)$ so that $g \cdot
\phi$ corresponds to the standard flag. Now let $\cM$ be a family of
$\cA_{t,l}$--modules.
By passing to a certain localization of $R$, we may assume $\cM$ is
an element in $\cR(\cF_{t,l}(\Delta); d)(R)$. We choose an element
$g$, so that $g \cdot \pi(\cM)$ is the standard flag $F(d)$ and
consider the family $g \cdot \cM$ equivalent to $\cM$. Using the
definition of the scheme $\cR(\cF_{t,l}(\Delta); d)$, the element
$\psi_{t-l-1}$ for $g \cdot \cM$ is an element in $\frp_u^{(l)}(d)
\otimes R$, the family we want to construct.

Note that a family over the ring of dual numbers is just a tangent
vector, both for $\frp_u^{(l)}(d)$ and $\cR(\cF_{t,l}(\Delta),d)$, and
two families are equivalent precisely if the corresponding tangent
vectors are in the same orbit under the group action. Both
constructions above preserve the group action and the equivalence classes of
the families.

Thus we obtain the first isomorphism below; the second isomorphism
is given by the Artin--Voigt lemma.
$$
\begin{array}{rcl}
T_x(\frp_u^{(l)}(d)) / T_x(P(d)\cdot x) & \cong &
T_{M(x)}(\cR(\cF_{t,l}(\Delta);d))/T_{M(x)}(\mGL(\Sigma d) \cdot (M(x))) \\
& \cong &
\mExt^1(M(x),M(x)).
\end{array}
$$
Here we consider the Zariski tangent space $T_x(\frp_u^{(l)}(d))$ of
$\frp_u^{(l)}(d)$ in $x$, the Zariski tangent space $T_x(P(d)\cdot
x)$ of the $P$--orbit of $x$ in $x$, the Zariski tangent space
$T_{M(x)}(\cR(\cF_{t,l}(\Delta);d))$ of $\cR(\cF_{t,l}(\Delta),d)$
in $M(x)$, and the Zariski tangent space of the orbit of $M(x)$.
This orbit consists precisely of the representations isomorphic to
$M(x)$.
\end{proof}

In particular, the proof shows that the scheme
$\cR(\cF_{t,l}(\Delta),d)$ is reduced. Let $\cM$ be a family of
modules over $\cA_{t,l}$ over $k[\varepsilon]/ \varepsilon^2$. If
the projection of the tangent direction of this family to the flag
variety  $\cR(\Aa_t; d)^{\mathrm{inj}}$ is not zero, we can find a
new equivalent family $\overline \cM$ having this property. Then
this gives a family of elements in $\frp_u^{(l)}(d)$, which we can
extend it to a family over $R = k[T]$, since $\frp_u^{(l)}(d)$ is an
affine space. Going back we obtain a family of  $\cA_{t,l}$--modules
equivalent to the original family $\cM$. So we can also extend $\cM$
to a family over $R = k[T]$ and $\cR(\cF_{t,l}(\Delta),d)$ is
reduced. This result is a little surprising, since the scheme of all
representations of $\cA_{t,l}$ is in general not reduced as follows
from \cite[Thm.\ 6.6]{schroerua}.
If the maps
$\phi_i$ are not assumed to be injective, we cannot prove such a
result.

Finally we note, that in the sequel, we only require one implication
of the theorem and apply it for $l=1$; that is we only use that the
codimension of $P \cdot x$ in $\frp_u^{(1)}$ is at most the
dimension of $\mExt^1(M(x),M(x))$. In our case, using a certain
filtration, it is possible to obtain upper bounds for the dimension
of the extension groups, whereas it seems difficult to obtain lower
bounds.


\section{Action with a dense orbit} \label{sectdense}

In this section we prove the direction of Theorem
\ref{maintheorem}(1) giving existence of a dense $B$--orbit in $\frn$ in
the stated cases.

Let $d$ be a thin dimension vector.  Let $B$ be the corresponding
Borel subgroup of $\mGL_n(k)$ and $\frn = \frp_u^{(1)}(d)$ the
corresponding ideal of $\frb$, see Notation \ref{notn}.  We
associate to $d$ the sequence $\underline a = \underline a (d)=
(a_0,a_1,\dots,a_r,a_{r+1})$ and set $e(d) = \sharp \{i \mid 1 \le i
\le r, \, a_i \text{ is even}\}$ as in Definition \ref{defna(d)}. We
prove that if $e(d)$ at most one, then $\frn$ is a prehomogeneous
space for $B$. In the next section we construct an explicit
representative $x$ of the dense orbit.

It is not feasible to prove that $B$ acts on $\frn$ with a dense
orbit using purely algebraic group theoretic results; such an
argument would be very technical for the case $e(d) = 1$.  Instead
we show that a certain $\Delta$-good module $M$ with
$\Delta$-dimension vector $d$ has no self-extensions. Consequently,
there is dense $B$--orbit in $\frn$ by Theorem \ref{Thmext}.

First we consider the case $e(d) = 0$.  In this case we decompose $d
= d^1 + d^2$ where
$$
d^1_i =
\left\{ \begin{array}{cl}   d_i & \text{if $i$ is odd} \\
                                0 & \text{if $i$ is even}
                                  \end{array} \right.
$$
and
$$
d^2_i =
\left\{ \begin{array}{cl}  d_i & \text{if $i$ is even} \\
                                  0 & \text{if $i$ is odd}
                                  \end{array} \right.
$$
Let $J$ and $K$ be the supports of $d^1$ and $d^2$ respectively. We
have corresponding standard modules $\Delta(J)$ and $\Delta(K)$ and
the condition $e(d) = 0$ implies that either $\Delta(J) \in
\{P(t-1),P(t)\}$ or $\Delta(K) \in \{P(t-1),P(t)\}$. Therefore, by
parts (2) and (3) of Theorem \ref{ThmDeltagoodproperties} we have
$\mExt^1(\Delta(J),\Delta(J)) = \mExt^1(\Delta(J),\Delta(K)) =
\mExt^1(\Delta(K),\Delta(J)) = \mExt^1(\Delta(K),\Delta(K)) = 0$.
Hence, setting $M = \Delta(J) \oplus \Delta(K)$ we have $\dim_\Delta
M = d$ and $\mExt^1(M,M) = 0$.  Therefore, $\frn$ is a
prehomogeneous space for $B$ by Theorem \ref{Thmext}.

\medskip

Now consider the case $e(d) = 1$, so $\underline a$ has one internal
even entry. The construction of a $\Delta$-good module $M$ with
$\dim_\Delta(M) = d$ and no self-extensions is more complicated in
this case. The idea is to construct $M$ as an extension of two
certain standard modules $S$ and $\R$. Using homological algebra we
show, that there exists a unique non-trivial extension
\begin{equation} \label{eqses}
0 \lra S \lra  M \lra \R \lra 0
\end{equation}
as an $\cA_{t,1}$--module and this module $M$ satisfies
$\mExt^1(M,M) = 0$ (Theorem \ref{Thmextvanishes}).

We start with the construction of the standard modules $\R =
\Delta(J)$ and $S = \Delta(K)$.  Let $b$ be the position of the end
of the internal even 1-string in $d$ and set
$$
J = \{\dots,b-5,b-3,b-1,b+2,b+4,b+6,\dots\}.
$$
Now define $K$ to be disjoint from $J$ and such that $J \cup K$ is
the support of $d$.  We let $c \in \{t-1,t\}$ be the greatest
element of $J$.

\begin{Exmp}
Let $t = 18$ and let $d = (1,1,0,1,0,1,0,1,1,1,0,1,1,1,1,0,1,1)$.
Then $b = 15$, $J = \{2,4,6,8,10,12,14,17\}$, $K = \{1,9,13,15,18\}$
and $c = 17$.
\end{Exmp}

We are now in a position to state the main theorem of this section.

\begin{Thm}\label{Thmextvanishes}
{\em (a)} $\mExt^1(\R,S) = k$ and $\mExt^1(S,\R) = 0$\\
{\em (b)} If $M$ is the unique non-trivial extension $0 \lra S \lra
M \lra \R \lra 0$, then $\mExt^1(M,M) = 0$.
\end{Thm}

This theorem along with Theorem \ref{Thmext} implies that $B$ acts
on $\frn$ with a dense orbit.  We prove the theorem in several steps
starting with an explicit projective resolution in Lemma
\ref{lemmaprojresolv}, and proceed with some vanishing results.

In the next lemma we use the notation for morphisms between standard
modules given in Definition \ref{defmor}.

\begin{Lemma}\label{lemmaprojresolv}
$0 \lra P(b) \overset{\psi}\lra P(b-1)\oplus P(c) \lra 0$ is a
projective resolution of $\R$, where $\psi = (\phi^{b-1},\phi^b)$.
\end{Lemma}

\begin{proof} Since $\R$ is a standard module, it is easy to obtain the
result using the covering of $\R$ and the covering of the projective
modules as explained in \cite{BHJAlg}.  Alternatively one can prove
the result using the explicit description of the projective modules
and the morphisms between them given in Section
\ref{sectpreliminary}.
\end{proof}

\begin{proof}[Proof of Theorem \ref{Thmextvanishes}] We first prove part
(a).  We have $\mExt^1(S,\R) = 0$ by Lemma \ref{lemdimext}. So we are
left to show that $\mExt^1(\R,S) = k$.

Applying $\mHom(-,S)$ to the projective resolution given in Lemma
\ref{lemmaprojresolv}, we obtain the long exact sequence
$$
0 \lra \mHom(\R,S) \lra \mHom(P(b-1),S) \oplus \mHom(P(c),S) \lra
\qquad \qquad \qquad ~
$$
$$
\qquad  \qquad \qquad \lra \mHom(P(b),S) \lra \mExt^1(\R,S) \lra 0.
$$
Let $l = \sharp(K \cap \{1,\dots,b-1\})$ and $m = \sharp (K \cap
\{1,\dots,c\})$. The dimensions of $\mHom$-groups given below can be
determined using Lemma \ref{lemdimhom}:
$$
\mhom(\R,S) = m, \, \mhom(P(b-1),S) = l, \, \mhom(P(c),S) = m,
$$
$$
\text{ and } \mhom(P(b),S) = l+1.
$$
Consequently, $\mext^1(\R,S) = 1$.

b) We show $\mExt^1(M,M) = 0$ in several steps using the exact
sequence \eqref{eqses} including $M$. First we note that by Theorem
\ref{ThmDeltagoodproperties}(3), we have $\mExt^1(\R,\R) = 0 =
\mExt^1(S,S)$.
\\
{\em Claim 1.} $\mExt^1(M,S) = 0$: \\
Applying $\mHom(-,S)$ to \eqref{eqses} we obtain a long exact
sequence containing
$$
\mHom(S,S) \lra \mExt^1(\R,S) \lra \mExt^1(M,S) \lra \mExt^1(S,S) = 0.
$$
Since the extension including $M$ is non-trivial, the first map in
the sequence above is not zero.  Therefore, it is surjective,
which implies that $\mExt^1(M,S) \cong \mExt^1(S,S) = 0$. \\
{\em Claim 2.} $\mExt^1(M,\R) = 0$: \\
Applying $\mExt^1(-,\R)$ to \eqref{eqses} we obtain an exact sequence
$$
0 = \mExt^1(\R,\R) \lra \mExt^1(M,\R) \lra \mExt^1(S,\R) = 0.
$$
{\em Claim 3.} $\mExt^1(M,M) = 0$: \\
Finally we have the exact sequence
$$
0= \mExt^1(M,S) \lra \mExt^1(M,M) \lra \mExt^1(M,\R) = 0,
$$
which is obtained by applying $\mHom(M,-)$ to \eqref{eqses}.  Hence,
$\mExt^1(M,M) = 0 $ as required.
\end{proof}

\begin{Rem}
Fix the standard subset $J$ of $\{1,\dots,t\}$ as above, let $K$ be
any standard subset.  Let $\R = \Delta(J)$ and $S = \Delta(K)$. We
note that the key properties of the standard subset $K$ for the
proof of Theorem \ref{Thmextvanishes} are $\mExt^1(S,\R) = 0$ and $b
\in K$.  There are additional possibilities for $S$ such that
$\mExt^1(S,\R) = 0$ provided by Lemma \ref{lemdimext}.  Therefore,
the same arguments can be used to prove a more general result that
$P(d)$ acts on $\frp_u^{(1)}(d)$ with a give a dense orbit, when $d$
may also contain some entries $2$.
\end{Rem}

\section{Representatives of the dense orbits}
\label{sectrepresentative}

In this section we construct a representative of the dense orbit of
$B$ in $\frn$ for the cases where such an orbit exists.  In Section
\ref{sectnotdense} we use these explicit constructions to prove that
a dense $B$--orbit does not exist in the cases given by Theorem
\ref{maintheorem}(1).  The construction is similar to that given in
\cite[\S 8]{BHRR} for Richardson elements in parabolic subgroups of
$\mGL_n(k)$.

Let $d$ be a thin dimension vector with at most one internal even
1-string.  Let $B$ be the corresponding Borel subgroup of
$\mGL_n(k)$ and $\frn$ the corresponding ideal of $\frb$, see
Notation \ref{notn}.  We also use the relabeled basis $\{f_j \mid j
\in \msupp d\}$ defined in Notation \ref{notn} and the elementary
matrices $f_{ij}$ with respect to this basis.

Below we construct two representatives of the dense $B$--orbit in
$\frn$. For the first representative $x$ it is easy to see the
module structure of the corresponding $\cA_{t,1}$--module $M(x)$. We
prove that the second representative $\bar x$ is minimal in the
sense of \cite[Defn.\ 7.1]{goo3}, this minimality in important for
the proofs in the subsequent sections.

We consider the case where $d$ has one internal even 1-string.  At
the end of this section (Remark \ref{remeasy}), we explain how to
simplify the construction in order to get a representative of the
dense $B$--orbit in the easier case when $d$ has no even internal
1-strings.

Let $b$ be the position of the end of the even 1-string in $d$ and
decompose the support of $d$ into standard subsets $J$ and $K$ as in
the previous section.

The construction of $x$ begins by defining a diagram $D(d)$ in the
plane.
\begin{enumerate}
\item Place vertices at the coordinates $(j,1)$ for $j \in J$ and
$(i,0)$ for $i \in K$.
\item From each vertex which is not leftmost in its row draw an
arrow to the vertex to its left.
\item Draw an arrow from the vertex at the point $(b+2,1)$ to the
vertex at $(b,0)$.
\item Label the vertices by their $x$-coordinates.
\item Let $x$ be the matrix which is the sum of the elementary
matrices $f_{ij}$ for which there is an arrow from vertex $j$ to
vertex $i$.
\end{enumerate}

We note that the Jordan normal form of the nilpotent matrix $x$ as
defined above consists of two Jordan blocks of size $\sharp J$ and
$\sharp K$.

We illustrate the construction of $x \in \frn$ in an example.

\begin{Exmp} \label{exD(d)}
Let $d = (1,1,0,1,0,1,1,0,1,1,1)$, so $\underline a = (2,1,2,3)$.
Then we have $b = 7$, $J = \{2,4,6,9,11\}$ and $K = \{1,7,10\}$ and we
construct the diagram
$$
\xymatrix@C=15pt{ & \bullet^2 & & \bullet^4 \ar[ll] & & \bullet^6 \ar[ll] &
& &
\bullet^9 \ar[lll] \ar[dll] & & \bullet^{11} \ar[ll] \\
\bullet^1 & & & & & & \bullet^7 \ar[llllll] & & & \bullet^{10}
\ar[lll]}
$$
This gives the representative
$$
x = \left(
\begin{array}{c}
  \includegraphics[width=3.5cm]{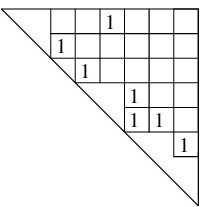}
\end{array}
\right).
$$
\end{Exmp}

To construct our second representative $\bar x$ we have to slightly
modify the diagram defining $x$ to obtain a diagram $\bar D(d)$. If
there is no vertex at a point $(i,0)$ with $i > b$, then no
modification to $D(d)$ is required. Otherwise let $i$ be minimal
such that $i > b$ and there is a vertex at $(i,0)$ in $D(d)$;
then we replace the arrow from vertex $i$ to vertex $b$ by an arrow from
$i$ to $b-1$.  Now $\bar x$ is defined from $\bar D(d)$ in
the same way that $x$ is defined from $D(d)$.

Now we show that $x$ and $\bar x$ are conjugate under $B$.  We may
assume that there is a vertex at a point $(i,0)$ with $i > b$, and
we let $i_1 < \dots < i_s$ be all such $i$.  By construction of $J$
and $K$ there must be at least $s$ vertices with $y$-coordinate $1$
that lie to the right of the vertex $b+2$.  Let $b+2 = j_1 < \dots <
j_s$ be the labels of the first $s$ such vertices. Consider the
change of base of $k^n$ given by $f_{i_l} \mapsto f_{j_l} -
f_{i_l}$, for $l = 1,\dots,s$.  It is clear that the matrix $g$
corresponding to this change of base lies in $B$ and that $gxg^{-1}
= \bar x$.

\begin{Exmp}
We continue with $d$ as in Example \ref{exD(d)}.  Then $\bar D(d)$
is as below
$$
\xymatrix@C=15pt{ & \bullet^2 & & \bullet^4 \ar[ll] & & \bullet^6 \ar[ll] &
& &
\bullet^9 \ar[lll] \ar[dll] & & \bullet^{11} \ar[ll] \\
\bullet^1 & & & & & & \bullet^7 \ar[llllll] & & & \bullet^{10}
\ar[ullll]}
$$
and determines the representative
$$
\bar x = \left(
\begin{array}{c}
  \includegraphics[width=3.5cm]{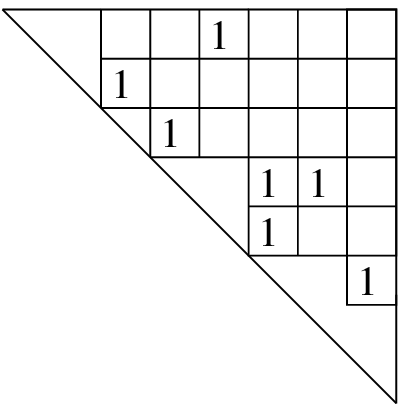}
\end{array}
\right).
$$
The matrix conjugating $x$ to $\bar x$ is
$$
g = \left(
\begin{array}{c}
  \includegraphics[width=3.5cm]{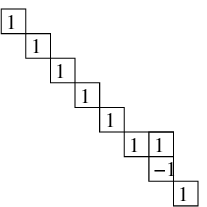}
\end{array}
\right).
$$
\end{Exmp}

Let $x,\bar x \in \frn$ be constructed from the diagrams $D(d)$ and
$\bar D(d)$ respectively.  In Proposition \ref{proprep} below we
prove that $x$ and $\bar x$ are representatives of the dense
$B$--orbit in $\frn$. Further, we prove that, with respect to a
certain enumeration of the positive roots of the root system of
$\mGL_n(k)$, we have that $\bar x$ is the minimal representative of
its $B$--orbit in the sense of \cite[Defn.\ 7.1]{goo3}, see also
\cite{verarr}.  In particular, this means that $\bar x$ (and $x$)
have minimal number of nonzero entries. In order to achieve this, we
need to give some notation and recall some terminology from
\cite{goo3}.  We warn the reader that in order to simplify our
proofs in the sequel some of the notation we give below is
nonstandard.  In particular, our labelling of the roots of
$\mGL_n(k)$ uses the relabeled basis $\{f_j \mid j \in \msupp d\}$
of $k^n$.

Let $T$ be the maximal torus of $\mGL_n(k)$ consisting of diagonal
matrices.  Let $\Phi = \{(i,j) \in \msupp d \times \msupp d
\mid i \ne j\}$ be the
root system of $\mGL_n(k)$ with respect to $T$, and let $\Phi^+ =
\{(i,j) \in \Phi \mid i < j\}$ be the system of positive roots
determined by $B$. For a root $\beta \in \Phi^+$ we write
$\frg_\beta$ for the corresponding root subspace of $\fru$. We
define $\Phi(\frn) = \{\beta \in \Phi^+ \mid \frg_\beta \subseteq
\frn\}$ to be the set of roots of $\frn$ with respect to $T$.

We define a linear order $\sqsubset$ on $\Phi^+$ by $(i,j) \sqsubset
(i',j')$ if $j < j'$, or $j = j'$ and $i > i'$. Then we enumerate
$\Phi(\frn) = \{\beta_1,\dots,\beta_{\dim \frn}\}$ so that $\beta_l
\sqsubset \beta_m$ if and only if $l < m$; we write $\beta_l =
(i_l,j_l)$. For example, if $\frn = \fru$ and $d = (1,1,0,1)$, then
we have $\beta_1 = (1,2), \beta_2 = (2,4), \beta_3 = (1,4), \beta_4
= (2,4)$. We define $B$--submodules $\frm_l$ of $\frn$ by $\frm_l =
\bigoplus_{m=l+1}^{\dim \frn} \frg_{\beta_m}$ for $l = 0,\dots,\dim
\frn$. Then we define $\frn_l = \frn/\frm_l$.

We study the orbits of $B$ in $\frn$ by considering the action of
$B$ on successive $\frn_l$.  Let $y \in \frn$ and let $\pi_l :
\frn_l \to \frn_{l-1}$ be the natural projection.  Consider the
fibre
$$
\pi_l^{-1}(y + \frn_{l-1}) =
\{y + \lambda e_{i_l,j_l} + \frm_l :
\lambda \in k \} \subseteq \frn_l.
$$
Then \cite[Lem.\ 5.1]{goo3} (see also \cite[Lem.\ 3.1]{verarr}) says
that for $y \in \frn$ either
\begin{enumerate}
\item[(I)] all elements of $\pi_l^{-1}(y + \frn_{l-1})$ are $U$--conjugate or
\item[(R)] no two elements of $\pi_l^{-1}(y + \frn_{l-1})$ are $U$--conjugate.
\end{enumerate}
We say that $l$ is an {\em inert} point of $y$ if (I) holds; and we
say that $l$ is a {\em ramification} point of $y$ if (R) holds.

In \cite[Defn.\ 7.1]{goo3} minimal representatives of $B$--orbits in
$\frn_l$ (with respect to an enumeration of $\Phi^+$) are defined.
We recall that \cite[Prop.\ 7.2 and Rem.\ 7.3]{goo3} says that a
$B$--orbit in $\frn_l$ contains a unique minimal representative.
Further, \cite[Lem.\ 5.5]{goo3} says that $y = \sum_{m = 1}^l a_m
e_{i_m,j_m} + \frm_l$ is the minimal representative of its
$B$--orbit in $\fru_l$ if and only if: $a_m = 0$ whenever $m$ is an
inert point of $y$; and if $a_m \ne 0$ and $\beta_m$ is linearly
independent of $\{\beta_{m'} \mid m' < m \text{ and } a_{m'} \ne
0\}$, then $a_m = 1$. From the above description of minimal
representatives of $B$--orbits in $\frn_l$, it follows that if $y +
\frm_l$ is the minimal representative of its $B$--orbit in $\frn_l$
and $l' < l$, then $y + \frm_{l'}$ is the minimal representative of
its $B$--orbit in $\frn_{l'}$.

The convention discussed in the following remark is required in the
proofs of both Propositions \ref{proprep} and \ref{Propnotdense}.

\begin{Rem} \label{remh_m}
Let $m \ne t$.  By definition of the function $\gamma$ in Notation
\ref{notn} $m$ is the $\gamma(m)$th lowest element of $\msupp d$. We
define $h_m \in \{1,\dots,M\}$ to be minimal so that $j_{h_m+1} >
m$. Let $B_{\gamma(m)}$ be the Borel subgroup of upper triangular
matrices in $\mGL_{\gamma(m)}(k)$, then we can identify $\frn_{h_m}$
with an ideal in the Lie algebra $\frb_{\gamma(m)}$ of
$B_{\gamma(m)}$. Under this identification we see that two elements
of $\frn_{h_m}$ are $B$--conjugate if and only if they are
$B_{\gamma(m)}$--conjugate. Also the identification allows us to
consider Jordan normal forms of elements of $\frn_{h_m}$.
\end{Rem}

We can now state the main proposition of this section.

\begin{Prop} \label{proprep}
Let $x,\bar x \in \frn$ be constructed from the diagrams $D(d)$ and
$\bar D(d)$ respectively.  Then $x$ and
$\bar x$ are representatives of the dense
$B$--orbit in $\frn$.  Moreover, $\bar x$ is the minimal
representative of this $B$--orbit.
\end{Prop}

\begin{proof}
Consider the $\cA_{t,1}$--module $M(x)$ corresponding to $x$ as
defined in Section \ref{sectdense}.  We see that the subspace of
$k^n$ spanned by $\{f_i : i \in K\}$ is stable under the action of
$x$.  It follows that there is a submodule $S$ of $M(x)$ isomorphic
to $\Delta(K)$. Further, we see that the quotient $\R = M(x)/S$ is
isomorphic to $\Delta(J)$. Now consider the direct sum of $\R \oplus
S$.  We can see $\R \oplus S = M(y)$, where $y = x - f_{b,b+2}$. Now
it is clear that $x$ is not $B$--conjugate to $y$, so $M(x)$ is not
isomorphic to $\R \oplus S$. Therefore, $x$ is a representative of
the dense $B$--orbit in $\frn$ by Theorems \ref{Thmext} and
\ref{Thmextvanishes}.  Since $\bar x$ is in the same $B$--orbit as
$x$ it follows that $\bar x$ is also a representative of the dense
$B$--orbit.

We define $A$ to be the subset of $\{1,\dots,\dim \frn\}$ such that
$\bar x = \sum_{l \in A} f_{i_l,j_l}$.  To show that $\bar x$ is the
minimal representative in its $B$--orbit it suffices to show that
each $l \in A$ is a ramification point of $\bar x$.

Let $l \in A$.  If $i_l = j_l - 2$, then it is clear that
$l$ is a ramification point of $x$, so assume that $l$ is such that
$i_l < j_l -2$.  Then the vertex $j_l$ has $y$-coordinate $0$,
or $j_l = b+2$ and $i_l = b-1$.  Now suppose for a
contradiction that $l$ is an inert point of $x$.  Then the minimal
representative of the $B$--orbit in $\frn_{h_{j_l}}$ is of the form
$$
y + \frm_{h_{j_l}} = x - f_{i_l,j_l} + \sum_{i \in \msupp d : i <
i_l} a_i f_{i,j_l} + \frm_{h_{j_l}},
$$
where $a_i \in \{0,1\}$.

We may identify $\frn_{h_{j_l}}$ with an ideal in
$\frb_{\gamma(j_l)}$ as discussed in Remark \ref{remh_m} and thus
view $x + \frm_{h_{j_l}}$ and $y + \frm_{h_{j_l}}$ as endomorphisms
of $k^{\gamma(j_l)}$.  Further, $x + \frm_{h_{j_l}}$ and $y +
\frm_{h_{j_l}}$ stabilize the subspace $k^{\gamma(i_l)-1}$ so induce
endomorphisms of $k^{\gamma(j_l)}/k^{\gamma(i_l)-1}$.
The kernels of these induced maps have different dimensions: the map
induced by $y+ \frm_{h_{j_l}}$ has kernel of dimension 3, whereas
the kernel of the map induced by $x + \frm_{h_{j_l}}$ is
2-dimensional. It follows that $x + \frm_{h_{j_l}}$ and $y +
\frm_{h_{j_l}}$ are not $B_{\gamma(j_l)}$--conjugate. Therefore, $x$
and $y$ cannot be $B$--conjugate. This contradiction implies that
$l$ is a ramification point of $x$, as required.
\end{proof}

\begin{Rem} \label{remeasy}
To end this section we mention how one can adapt the above
construction to give representatives of the dense $B$--orbit in
$\frn$ in case $d$ has no internal even 1-strings.  In this case we
decompose the support of $d$ in to subsets $J$ and $K$ as in the
previous section.  We construct a diagram using the same method as
for $D(d)$ except that we omit step 3, then use this diagram to
define $x \in \frn$. One can prove that $x$ is a representative in
the dense $B$--orbit in $\frn$ and that it is the minimal
representative of this orbit using arguments as in the proof of
Proposition \ref{proprep}.
\end{Rem}

\section{Actions without dense orbit}\label{sectnotdense}

In this section we prove non-existence of a dense $B$--orbit in the
cases given by Theorem \ref{maintheorem}(1). The proof requires the
explicit constructions of representatives of dense $B$--orbits (when
they exist) given in Section \ref{sectrepresentative}.

Let $d$ be a thin dimension vector.  Let $B$ be the corresponding
Borel subgroup of $\mGL_n(k)$ and $\frn$ the corresponding Lie ideal
of $\frb$, see Notation \ref{notn}.  We define the sequence
$\underline a = \underline a(d)$ and the number $e(d)$ as in
Definition \ref{defna(d)}.  For the proof of the following
proposition we use the terminology for $B$--orbits given in the
previous section.  In particular, we use the relabeled standard
basis $\{f_j \mid j \in \msupp d\}$ of $k^n$ and we enumerate the
roots $\Phi(\frn)$ of $\frn$ as before. Further, we define the
quotients $\frn_l$ of $\frn$ and we have minimal representatives of
$B$--orbits in $\frn_l$.

\begin{Prop}\label{Propnotdense}
Suppose $e(d) \ge 2$. Then $B$ does not act on $\frn$ with a dense
orbit.
\end{Prop}

\begin{proof} Let $\frm$ be an ideal of $\frb$ contained in $\frn$.
If $B$ acts on $\frn$ with a dense orbit, then $B$ also acts on
$\frn/\frm$ with a dense orbit.  It follows that it suffices to
consider the case where we have $\underline a =
(1,a_1,\ldots,a_r,1)$ with $a_1$ and $a_r$ even and $a_i$ odd for $2
\le i \le r-1$.

We have $\gamma(t-2) = n-1$.  Define $h_{t-2}$ and $B_{\gamma(t-2)}
= B_{n-1}$ as in Remark \ref{remh_m}. We can identify
$\frn_{h_{t-2}}$ with an ideal in $\frb_{n-1}$, it corresponds to
the dimension vector $\hat d = (d_1,\dots,d_{t-1})$ and the 1-string
sequence $\hat a = (1,a_1,\ldots,a_r)$.  We know there exists a
dense orbit for the action of $B$ on $\frn_{h_{t-2}}$ and this dense
orbit has minimal representative $\bar x + \frm_{h_{t-2}} = \sum_{l
\in A} f_{i_l,j_l} + \frm_{h_{t-2}}$ that can be constructed using
the diagram $\bar D(\hat d)$ as described in the previous section.
The Jordan normal form of $\bar x + \frm_{h_{n-1}}$ consists of two
Jordan blocks of size $\lambda \ge \mu$, say.  It is easy to see
from the construction of $\bar x + \frn_{h_{n-1}}$ that we must in
fact have $\lambda > \mu$.

Suppose there is a dense $B$--orbit in $\frn$.  The minimal
representative of this orbit is of the form $\bar x = \sum_{l \in A}
f_{i_l,j_l} + \sum_{i=1}^{n-1} c_i e_{i,n}$, where $c_i \in
\{0,1\}$. It is clear that we must have that $c_{n-1} = 1$, i.e.\
$h_{n-1}+1$ is a ramification point of $\bar x$.  It now follows
from \cite[Lem.\ 5.7 and Prop.\ 7.7]{goo3} that, for $\bar x$ to be
a representative of a dense $B$--orbit in $\frn$, we must have that
$l$ is an inert point of $x$ for $l = h_{n-1}+2,\dots,M$; so that
$c_i = 0$ for $i =1,\dots,n-2$. However, we see that the Jordan
normal form of $\bar x$ has blocks of size $\lambda \ge \mu + 1$,
and the Jordan decomposition of $\bar x + e_{n-2,n}$ is given by the
partition $\lambda + 1 \ge \mu$. The closure ordering for conjugacy
classes of nilpotent matrices is well-known (see for example
\cite[Thm.\ 3.10]{hess}) and this ordering tells us that $x +
e_{n-2,n}$ is not contained in the closure of the $\mGL_n(k)$--orbit
of $\bar x$ and so certainly cannot be contained in the closure of
the $B$--orbit of $x$.  This contradiction implies that there cannot
be a dense $B$--orbit in $\frn$.
\end{proof}

\section{The codimension of a maximal $B$--orbit}\label{sectparameters}

In this section we improve our results and compute the codimension
of a $B$--orbit of maximal dimension in the ideal $\frn$.  In fact
we construct a dense family of $B$--orbits in $\frn$ of maximal
dimension.

Let $d$ be a thin dimension vector.  Let $B$ be the corresponding
Borel subgroup of $\mGL_n(k)$ and $\frn = \frp_u^{(1)}(d)$ the
corresponding Lie ideal of $\frb$, see Notation \ref{notn}.  We
define the number $e(d)$ as in Definition \ref{defna(d)}. We denote
by $\mcd(B;\frn)$ the codimension in $\frn$ of a $B$--orbit of
maximal dimension.

The main result in this section is the following theorem, which
gives part (2) of Theorem \ref{maintheorem}.

\begin{Thm}\label{Thmcodimmax}
Let $d$ be a thin dimension vector.  Then
$$
\mcd(B;\frn) = \max\{0,e(d)-1\}.
$$
\end{Thm}

For the remainder of this section we assume that $e(d) \ge 2$.  The
cases $e(d) \le 1$ in Theorem \ref{Thmcodimmax} follow from the
results of Section \ref{sectdense}.  We prove the theorem in two
steps. Using homological algebra, we compute an upper bound for
$\mcd(B;\frn)$. That is we show that for certain
$\cA_{t,1}$--modules, we have $\mExt^1(M,M) \le e(d) - 1$.  We
complete the proof by constructing an $(e(d) - 1)$-dimensional
affine family such that the geometric fibres $\bar X$ of the family
consists of pairwise non-conjugate elements $x \in \frn$, such that
for each $x \in \bar X$ the corresponding module $M(x)$ is
isomorphic to a module $M$ for which we have shown $\mExt^1(M,M) \le
e(d) - 1$.

\medskip

We obtain the upper bound $\mcd(B;\frn) \le e(d)-1$ from Lemma
\ref{extlemma} below and Theorem \ref{Thmext}.  First we introduce
some notation that is required for the proof of this lemma.

Let $b_1 < \dots < b_{e(d)}$ be the ends of the internal even
1-strings in $d$.  Define $J =
\{\dots,b_1-1,b_1+2,\dots,b_2-1,b_2+2,\dots,b_3-1,b_3+2,\dots\}$ to
be the largest standard subset containing $b_i-1,b_i+2$ and all
possible $b_1-1-2j$ and $b_{e(d)}+2 +2j$ for $j\geq 0$. Further, let
$c$  be the greatest element of $J$ and let $K$ be the subset of
$\{1,\dots,t\}$ such that $\msupp d$ is the disjoint union of $J$
and $K$. Define the standard modules $\R = \Delta(J)$ and $S =
\Delta(K)$.

\begin{Exmp} \label{exmpJK}
Let $d = (1,0,1,1,0,1,1,1,1,0,1,1,1,0,1)$.  Then we have $J =
\{1,3,6,8,11,13,15\}$ and $K = \{4,7,9,12\}$.
\end{Exmp}

\begin{Lemma}\label{extlemma}
{\em (a)} $\mext^1(\R,S) = e(d)$\\
{\em (b)} $\mext^1(S,\R) = 0$ \\
{\em (c)} For any non-trivial extension $M$ of $\R$ and $S$ we have
$\mext^1(M,M) \le e(d)-1$.
\end{Lemma}

\begin{proof} (a)
We first consider a projective resolution of $\R$
$$
0 \lra P(b_1) \oplus \ldots \oplus P(b_m) \overset{\psi} \lra P(b_1
- 1) \oplus \ldots \oplus P(b_{e(d)} - 1) \oplus P(c),
$$
where $\psi = \psi_1 \oplus \dots \psi_{e(d)}$, with $\psi_i =
(0,\dots,0,\phi^{b_i - 1},\phi^{b_i},0,\dots,0)$ for $i =
1,\dots,e(d)-1$ and $\psi_{e(d)} = (0,\dots,0,\phi^{b_{e(d)} -
1},\phi^{b_{e(d)}})$, here we use the notation of Definition
\ref{defmor}. As in the proof of Theorem \ref{Thmextvanishes}, this
gives rise to a long exact sequence.
$$
0 \lra \mHom(\R,S) \lra \mHom(P(b_1-1),S) \oplus \dots
\oplus \mHom(P(c),S) \lra \qquad \qquad \qquad ~
$$
$$
\qquad  \qquad \qquad \lra \mHom(P(b_1),S) \oplus \dots \oplus
\mHom(P(b_{e(d)}),S) \lra \mExt^1(\R,S) \lra 0.
$$
Computing dimensions of $\mHom$-groups using Lemma \ref{lemdimhom}
yields
$$
\mhom(\R,S) = \mhom(P(c),S), \text{ and }
$$
$$
\mhom(P(b_i),S) = \mhom(P(b_i-1),S) + 1 \text{ for $i = 1,\dots,e(d)$.}
$$
From this we can calculate $\mext^1(\R,S) = e(d)$.

(b) Using Lemma \ref{lemdimext}, we get $\mExt^1(S,\R) = 0$.

(c) We consider a non-trivial extension
\begin{equation} \label{eqses2}
0 \lra S \lra M \lra \R \lra 0.
\end{equation}
Applying $\mHom(-,S)$ to the short exact sequence above we obtain a
long exact sequence containing
$$
\mHom(S,S) {\lra} \mExt^1(\R,S) \lra \mExt^1(M,S) \lra 0,
$$
since $\mExt^1(S,S) = 0$ by Theorem \ref{ThmDeltagoodproperties}(4).
Moreover, since the first map in the exact sequence above is
not the zero map, we obtain $\mext^1(M,S) \le e(d)-1$. Applying
$\mHom(-,\R)$ to \eqref{eqses2} we obtain the exact sequence
$$
0 = \mExt^1(\R,\R) \lra \mExt^1(M,\R) \lra \mExt^1(S,\R) = 0
$$
and deduce that $\mExt^1(M,\R) = 0$. Finally applying $\mHom(M,-)$ to
\eqref{eqses2} we obtain the exact sequence
$$
\mExt^1(M,S) \lra \mExt^1(M,M) \lra \mExt^1(M,\R) = 0
$$
implying that $\mExt^1(M,M) \le e(d)-1$ as required.
\end{proof}

We now explain a construction of a family of elements $x \in \frn$,
for which $M(x)$ is a nontrivial extension of $S$ by $\R$.  This
construction generalizes the one given in Section
\ref{sectrepresentative} and begins by defining a diagram in the
plane. We use the relabeled basis $\{f_j \mid j \in \msupp d \}$ of
$k^n$ from Notation \ref{notn}. We denote the diagram in the plane
in this construction by $D(d)$.

\begin{enumerate}
\item Place vertices at the coordinates $(j,1)$ for $j \in J$ and
$(i,0)$ for $i \in K$.
\item From each vertex which is not leftmost in its row and is not at
$(b_i+2,1)$ for some $i = 2,\dots,e(d)$, draw an arrow to the vertex
on its left.
\item For $i = 1,\dots,e(d)$, draw an arrow from the vertex at the point
$(b_i+2,1)$ to the vertex at the point $(b_i,0)$.
\item For $i = 2,\dots,e(d)$ draw a twiddly arrow from the vertex at
$(b_i+2,1)$ to the vertex at $(b_i-1,1)$ and label this twiddly
arrow $i$.
\item Number the vertices by their $x$-coordinate.
\item Let $A \subset \{(i,j) \mid i < j\}$ be the set of all pairs
such that there is a non-twiddly arrow from $j$ to $i$.
\item Define $\cF = \sum_{(i,j) \in A} f_{ij} + \sum_{i=2}^{e(d)} x_i
f_{b_i-1,b_i+2} \in \frn \otimes k[x_2,x_2^{-1},\dots,x_m,
x_m^{-1}]$, where $x_2,\dots,x_m$ are indeterminates. Then $\cF$ is
a family of elements of $\frn$ over the affine algebra
$k[x_2,x_2^{-1},\dots,x_m, x_m^{-1}]$.
\end{enumerate}

\begin{Exmp} \label{exmpD(d)}
Let $d = (1,0,1,1,0,1,1,1,1,0,1,1,1,0,1)$ as in Example
\ref{exmpJK}.  Then the diagram $D(d)$ is given below
$$
\xymatrix@C=9pt{\bullet^1 & & \bullet^3 \ar[ll] & & & \bullet^6
\ar[lll] \ar[dll] & & \bullet^8 \ar[ll] & & & \bullet^{11}
\ar@{~>}[lll]_2 \ar[dll] &
& \bullet^{13} \ar[ll] & & \bullet^{15} \ar[ll] \\
 & & & \bullet^4 & & & \bullet^7 \ar[lll] & & \bullet^9 \ar[ll] & &
& \bullet^{12}
 \ar[lll] }
$$
We have
$$
x = f_{1,3} + f_{3,6} + f_{4,6} + f_{4,7} + f_{6,8} + f_{7,9} +
x_2 f_{8,11} + f_{9,11} + f_{9,12} + f_{11,13} + f_{13,15}.
$$
\end{Exmp}

Each $(t_2,\dots,t_m) \in (k^\times)^m$ defines a maximal ideal of
$k[x_2,x_2^{-1},\dots,x_m,x_m^{-1}]$, and we write
$\cF(t_2,\dots,t_m) = \sum_{(i,j) \in A} f_{ij} + \sum_{i=2}^m t_i
f_{b_i-1,b_i+2} \in \frn$.

The following lemma is proved in the same way as the first part of Proposition
\ref{proprep}, so we omit the proof.

\begin{Lemma} \label{lemmafamily}
Let $\R$, $S$ and $\cF$ be as above and let $(t_2,\dots,t_m) \in
(k^\times)^m$.  Then $M(\cF(t_2,\dots,t_m))$ is a nontrivial
extension of $S$ by $\R$.
\end{Lemma}

As in Section \ref{sectrepresentative} we require an alternative
construction in order to obtain minimal representatives of
$B$--orbits. The modification required is explained below.

To construct our second family $\bar \cF$ we have to slightly modify
the diagram $D(d)$ to obtain a diagram $\bar D(d)$. For each $i = 1,
\dots,e(d)-1$, we let $j_i$ be the label of the vertex to the right
of $b_i$.  Replace the arrow from $j_i$ to $b_i$ by an arrow from
$j_i$ to $b_i-1$. If there is a vertex $j_{e(d)}$ to the right of
the vertex $b_{e(d)}$, then we replace the arrow from $j_{e(d)}$ to
$b_{e(d)}$ by an arrow from $j_{e(d)}$ to $b_{e(d)}-1$.  Now the
family $\bar \cF \in \frn \otimes
k[x_2,x_2^{-1},\dots,x_m,x_m^{-1}]$ is defined from $\bar D(d)$ in
the same way that $\cF$ is defined from $D(d)$.

Let $(t_2,\dots,t_{e(d)}) \in (k^\times)^{e(d)-1}$.  Then it is
straightforward to show that $\cF(t_2,\dots,t_{e(d)})$ is
$B$--conjugate to $\bar \cF(t_2,\dots,t_{e(d)})$.  This can be done
inductively using the arguments from Section
\ref{sectrepresentative} showing that, in the situation of that
section, the alternative representatives $x$ and $\bar x$ of the
dense $B$--orbit are conjugate.

\begin{Exmp}
Let $d$ as in Examples
\ref{exmpJK} and \ref{exmpD(d)}.  Then the diagram $\bar D(d)$ is
$$
\xymatrix@C=9pt{\bullet^1 & & \bullet^3 \ar[ll] & & & \bullet^6
\ar[lll] \ar[dll] & & \bullet^8 \ar[ll] & & & \bullet^{11}
\ar@{~>}[lll]_2 \ar[dll] &
& \bullet^{13} \ar[ll] & & \bullet^{15} \ar[ll] \\
 & & & \bullet^4 & & & \bullet^7 \ar[ullll] & & \bullet^9 \ar[ll] & &
& \bullet^{12}
 \ar[ullll] }
$$
We have
$$
\bar \cF = f_{1,3} + f_{3,6} + f_{3,7} + f_{4,6} + f_{6,8} + f_{7,9}
+ x_2 f_{8,11} + f_{8,12}  + f_{9,11} + f_{11,13} + f_{13,15}.
$$
\end{Exmp}

For the statement Lemma \ref{lemfammin}, we use the terminology of
minimal representatives with the enumeration of $\Phi(\frn)$ given
in Section \ref{sectrepresentative}.  The lemma is proved using
arguments from the proof of Proposition \ref{proprep}, so we omit
its proof.

\begin{Lemma} \label{lemfammin}
For each $(t_2,\dots,t_{e(d)}) \in (k^{\times})^{e(d)-1}$ have that
$\bar \cF(t_2,\dots,t_{e(d)})$ is the minimal representative of its
$B$--orbit.
\end{Lemma}

We are now in a position to prove Theorem \ref{Thmcodimmax}.

\begin{proof}[Proof of Theorem \ref{Thmcodimmax}]
Let $\bar X = \{\bar \cF(t_2,\dots,t_m) \in \frn \mid
(t_2,\dots,t_m) \in (k^\times)^{e(d)-1}\}$ be the geometric fibres
of the family $\bar \cF$.  The elements of $\bar X$ are minimal
representatives in their $B$--orbits by Lemma \ref{lemfammin}.
Therefore, the uniqueness of minimal representatives (\cite[Prop.\
7.2 and Rem.\ 7.3]{goo3}) means that no two elements of $\bar X$ lie
in the same $B$--orbit.  By Lemmas \ref{extlemma} and
\ref{lemmafamily}, and Theorem \ref{Thmext}, the codimension in
$\frn$ of the $B$--orbit of any $y \in \frn$ is $e(d)-1$ or less. It
follows that $B \cdot \bar X = \bigcup_{y \in \bar X} B \cdot y$ has
dimension $\dim \frn$ and is therefore dense in $\frn$. It follows
from \cite[Cor.\ AG.10.3]{bor} that $Z = \{z \in \frn \mid
\mcodim_\frn B \cdot z < e(d)-1\}$ is open in $\frn$. If $Z$ is
nonempty, then its intersection with $B \cdot \bar X$ is nonempty.
Then $Z \cap \bar X$ is open and nonempty in $\bar X$. But then the
dimension of $B \cdot (Z \cap \bar X)$ is greater than the dimension
of $\frn$, which is absurd. Hence, $Z$ is empty and we have
$\mcd(B;\frn) = e(d)-1$.
\end{proof}

Next we give an example of a dense family of $B$--orbits.

\begin{Exmp} For $d=(1,0,1,1,1,0,1,0,1,0,1,0,1,1,0,1,0,1,1,1,1,0,1,1)$ we
have a family with one parameter $t$.  This family is illustrated in
Figure \ref{figfamily} below.

\begin{figure}[h!tb]

\begin{center}

\includegraphics[height=5cm]{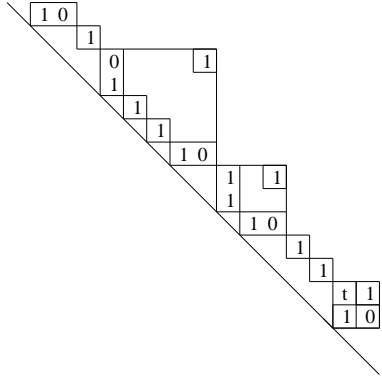}

\caption{A family of $B$--orbits} \label{figfamily}

\end{center}

\end{figure}

\end{Exmp}

\section{Maximal $B(q)$--conjugacy classes in $N(q)$}

We assume, for this section, that $k$ is the algebraic closure of
the finite field of $p$ elements $\Ff_p$, where $p$ is a prime.  Let
$d$ be a thin dimension vector.  Let $B$ and $\frn$ be the be the
corresponding Borel subgroup of $\mGL_n(k)$ and Lie ideal of $\frb$,
see Notation \ref{notn}. We write $N$ for the unipotent normal
subgroup of $B$ with Lie algebra $\frn$.

We note that $B$ has a standard definition over $\Ff_p$ and $N$ is a
subgroup of $B$ defined over $\Ff_p$. For a power $q$ of $p$ we
write $B(q)$ for the group of $\Ff_q$-rational points in $B$, that
is $B(q) = B \cap \mGL_n(q)$; likewise we write $N(q)$ and $\frn(q)$
for the $\Ff_q$-rational points of $N$ and $\frn$ respectively. The
purpose of this final section is to deduce a result (Theorem
\ref{thmfinite}) about the maximal size of a $B(q)$-conjugacy class
in $N(q)$ from Theorem \ref{Thmcodimmax}. This is achieved using
results of the first author from \cite{goo3}.

We associate to $d$ the positive integer $e(d)$ as in Definition
\ref{defna(d)}.  We use the terminology for $B$--orbits recalled
from \cite{goo3} in Section \ref{sectrepresentative}.  In
particular, we enumerate $\Phi(\frn)$ as in Section
\ref{sectrepresentative}, then we have minimal representatives of
$B$--orbits in $\frn$ and the notion of inert and ramifications
points.

Let $y = \sum_{(i,j) \in A} a_i f_{i,j} \in \frn(q)$ be the minimal
representative of its $B$--orbit in $\frn$, where $A$ is a subset of
$\Phi(\frn)$ and $a_i \in \Ff_q^\times$.  We write $\mIn(y)$ for the
number of inert points of $y$.  By \cite[Prop.\ 7.7]{goo3} we have a
factorization $C_B(y) = C_T(y)C_U(y)$ and this induces a
factorization $C_{B(q)}(y) = C_{T(q)}(y)C_{U(q)}(y)$.  We have that
$\dim U \cdot y = \mIn(y)$ by \cite[Lem.\ 5.3]{goo3} and we have
$|U(q) \cdot y| =  q^{\mIn(y)}$ by \cite[Prop.\ 6.4]{goo3}. Let
$\mm(y)$ be the rank of the sublattice of the root lattice of $\Phi$
generated by $A$.  Then it is easy to see that $\dim T \cdot y =
\mm(y)$ and $|T(q) \cdot y| = (q-1)^{\mm(y)}$.  Putting this all
together we have that $\dim B \cdot y = \mm(y) + \mIn(y)$ and $|B(q)
\cdot y| = (q-1)^{\mm(y)}q^{\mIn(y)}$.

We are now in a position to prove the theorem of this section.

\begin{Thm} \label{thmfinite}
The maximal size of a $B(q)$--conjugacy class in $N(q)$ is: \\
{\em (a)} $(q-1)^{n-2} q^{\dim N - (n-2)}$, if $e(d) = 0$; \\
{\em (b)} $(q-1)^{n-1} q^{\dim N - (n - 2 + e(d))}$, if $e(d) > 0$.
\end{Thm}

\begin{proof}
We first note that the map $x \mapsto 1 + x$ from $\frn$ to $N$ is a
$B$-equivariant and sends $\frn(q)$ to $N(q)$.  Therefore, we can
consider the action of $B$ on $\frn$ rather than on $\frn$.

We first consider the case $e(d) = 1$.  Let $\bar x$ be the minimal
representative of the dense $B$--orbit in $\frn$ constructed in
Section \ref{sectrepresentative}.  The entries of $\bar x$ are all 0
or 1 so we have $\bar x \in \frn(q)$.   Since $B \cdot \bar x$ is
dense in $\frn$, we have $\dim B \cdot \bar x = \dim \frn$.  It is
easy to see that $\mm(\bar x) = n-1$ and it follows that $\mIn(\bar
x) = \dim \frn - (n-1)$.  Therefore, we have $|B(q) \cdot \bar x| =
(q-1)^{n-1} q^{\dim N - (n - 1)}$.

Now $B \cdot \bar x$ is open and dense $\frn$ and the $U$--orbit of
any $z \in B \cdot \bar x$ has dimension $\dim \frn - (n-1)$.  It
follows from \cite[Cor.\ AG.10.3]{bor} that $\dim U \cdot z \le \dim
\frn - (n-1)$ for all $z \in \frn$.  It is clear that $\dim T \cdot
z \le n-1$ for any $z \in \frn$.  Therefore, if $y \in \frn(q)$ is
the minimal representative of its $B$--orbit in $\frn$ we have
$\mm(y) \le n-1$ and $\mIn(y) \le \frn - (n-1)$, so that $|B(q)
\cdot y| \le (q-1)^{n-1} q^{\dim \frn - (n - 1)}$.  The case $e(d) =
1$ in the theorem now follows, because each $B$--orbit in $\frn$ has
a minimal representative.

The case $e(d) = 0$ can be proved in more or less the same way.  Let
$x \in \frn(q)$ be the minimal representative of the dense
$B$--orbit as constructed in Section \ref{sectrepresentative}. Then
we can show that $|B(q) \cdot x| = (q-1)^{n-2} q^{\dim N - (n-2)}$.
In this case it is necessary, to observe that $\frn \setminus (B
\cdot x)$ has dimension strictly less than $\dim \frn$. It follows
that if $y \in \frn \setminus (B \cdot x)$ is the minimal
representative of its $B$-orbit, then $\dim B \cdot y = \mm(y) +
\mIn(y) < \dim \frn$, which forces $|B(q) \cdot y| < (q-1)^{n-2}
q^{\dim N - (n-2)}$ (an analogous strict inequality could have also
been shown in the case $e(d) = 1$).

The case $e(d) > 1$ can be proved in much the same way as for $e(d)
= 1$. Let $\bar \cF \in \frn \otimes
k[x_2,x_2^{-1},\dots,x_{e(d)},x_{e(d)}^{-1}]$ be the family of
elements of $\frn$ from Section \ref{sectparameters}. Then for
$(t_2,\dots,t_{e(d)}) \in \Ff_q^{\times}$, we have that
$$
|B(q) \cdot \bar \cF(t_2,\dots,t_{e(d)})| = (q-1)^{n-1} q^{\dim \frn
- (n - 2 + e(d))}.
$$
Since, $Z = \bigcup_{(t_2,\dots,t_{e(d)}) \in k^{\times}} B \cdot
\bar \cF(t_2,\dots,t_{e(d)})$ is dense in $\frn$ and the $U$--orbit
of each element of $Z$ has dimension $\dim \frn - (n - 2 + e(d))$,
we can deduce as for the case $e(d) = 1$ that $(q-1)^{n-1} q^{\dim
\frn - (n - 2 + e(d))}$ is an upper bound for the size of a
$B(q)$--orbit in $\frn(q)$.
\end{proof}

\medskip

{\small \noindent
Simon M.~Goodwin\\
New College \\ Oxford \\ OX1 3LB
\\ Email: simon.goodwin@new.ox.ac.uk \\ http://home.imf.au.dk/goodwin/}

\medskip

{\small \noindent
Lutz Hille\\
Mathematisches Seminar\\
Universit\"at Hamburg\\
D-20 146 Hamburg\\
Email: hille@math.uni-hamburg.de\\
http://www.math.uni-hamburg.de/home/hille/}

\end{document}